\documentclass[twocolumn,10pt]{article}

\usepackage[a4paper,top=25mm,bottom=25mm,left=18mm,right=18mm]{geometry}
\usepackage{graphicx}
\usepackage{amsmath,amssymb}
\usepackage{booktabs}
\usepackage{multirow}
\usepackage{url}

\newcommand{\ebar}{\bar{e}}

\newcommand{\epsa}{\varepsilon_{\mathrm{a}}}
\newcommand{\epsr}{\varepsilon_{\mathrm{r}}}

\title{Performance Evaluation of an Adaptive Quadrature and a Double Exponential Formula Using Arbitrary-Precision Floating-Point Arithmetic}
\author{Tomonori Kouya\thanks{Faculty of Science and Engineering, Otemon Gakuin University}}
\date{}

\begin{document}

\maketitle

\begin{abstract}
Using arbitrary-precision arithmetic provided by the GNU Multiple Precision Floating-Point Reliable Library, we implement, for an arbitrary mantissa length, two existing methods: the doubly adaptive AQE11D of Hibino et al.\ and Takahasi and Mori's double exponential (DE) formula. We evaluate them for Kahaner's 21 test problems. For both absolute tolerances $10^{-50}$ and $10^{-100}$, AQE11D attains target accuracy on all 21 problems; however, for strong endpoint singularity such as $1/\sqrt{x}$, it requires about $5.4\times10^{7}$ function evaluations at $10^{-100}$, roughly $7\times10^{4}$ times as many as the DE formula. The formula converges on 18 problems at both tolerances, demonstrating its strength against endpoint singularities but also its failure, as it stands, on problems with a singularity inside the integration interval.
\end{abstract}


\section{Introduction}
Ill-conditioned problems, the evaluation of special functions, and numerical verification all call for numerical integration beyond IEEE 754 binary64 precision. The arbitrary-precision floating-point library GNU MPFR~\cite{mpfr} is built on the arbitrary-precision natural-number kernel
(MPN) of GMP~\cite{gmp} and provides correctly rounded arithmetic operations, elementary functions, and special functions; it is the de facto standard library for arbitrary-precision real arithmetic.

In this study, we implement the adaptive quadrature AQE11D and the double
exponential (DE) formula on a CPU using MPFR in order to compare them on Kahaner's 21 test problems~\cite{kahaner1971} in attained accuracy, number of function evaluations, serial execution time, and OpenMP parallel performance. We thereby propose design guidelines for a numerical integrator that meets a user-specified accuracy requirement.

Ninomiya~\cite{ninomiya1980} added two points to the 9-point Newton--Cotes rule and constructed an error estimation functional $\ebar$ from the resulting 11 function values; from the way $\ebar$ changes when the
interval width is halved, his adaptive routines AQNN5D/AQNN7D/AQNN9D estimate discontinuities, logarithmic singularities, and algebraic singularities. Hibino et al.~\cite{hibino2003,hasegawa2007} combined this 11-point scheme with
a sequence of higher-order rules of Favati--Lotti--Romani type~\cite{flr1991} to obtain the doubly adaptive AQE11D, the existing method implemented here; neither method is proposed in this note. By contrast, Takahasi and
Mori's~\cite{takahasi1974} DE formula is a representative method that converges rapidly for integrals with endpoint singularities. Kahaner's 21 test problems are widely used for comparing one-dimensional quadrature methods and in Gonnet's~\cite{gonnet2010} and Espelid's~\cite{espelid2003} evaluations of adaptive quadrature.

In our implementation, the coefficients and abscissae of the quadrature rules are constructed as MPFR variables. Since the coefficients of~\cite{ninomiya1980,flr1991} are rational numbers and their abscissae are binary fractions, this method is highly compatible with arbitrary-precision binary floating-point arithmetic. The rules and their error control are as published in~\cite{ninomiya1980,hibino2003,hasegawa2007,flr1991}; what is new here is marked in Sects.~2.1--2.2.

\section{Arbitrary-Precision Arithmetic and the Implemented Algorithms}
Here $\epsa$ and $\epsr$ are the absolute and relative error tolerances requested by the user for $I=\int_a^b f\,dx$, $P_B$ is the mantissa length in bits, and Nfe the number of evaluations of $f$. All floating-point operations herein are represented by MPFR's \texttt{mpfr\_t} with the mantissa precision fixed at $P_B$ bits, and the rounding mode is round-to-nearest (\texttt{MPFR\_RNDN}). The MPFR
API (\texttt{mpfr\_add}, \texttt{mpfr\_mul}, \texttt{mpfr\_div}, \texttt{mpfr\_exp}, etc.)\ is called directly through a C++ value-type wrapper. The rational coefficients of the quadrature rules and the
abscissae of the form $j/16$ are generated with ordinary integer arithmetic together with \texttt{mpfr\_mul\_2si}. This way, the values can be constructed to the required precision without going through decimal string representations of the coefficients.

The errors are evaluated with 800-bit MPFR arithmetic against reference values of 210 digits: analytic for the 15 problems with an elementary antiderivative (among them problem 7, $\int_0^1 x^{-1/2}dx=2$), and mpmath~\cite{mpmath} at 240 digits, checked at 400, for the other six (5, 8, 12, 13, 17, 18). We describe the algorithms of AQE11D and the DE formula in the following sections.

\subsection{Adaptive Quadrature AQE11D}
Consider an interval $[a,a+2h]$ of half-width $h$. A 9-point rule $Q$ is constructed from the function values at the nine equally spaced points and the error estimation functional $\ebar$ from the 11 function values that are derived after adding two points located at 16ths of the interval. The basic procedure includes the following: evaluating the initial 11 points, bisecting the interval, evaluating the six new points, testing for
convergence, testing for an anomalous point, and applying semi-analytic processing when necessary. The convergence test is as follows:

\begin{equation}
 |e| \leq
 \max\!\left(\epsa,\epsr|S'|\right)
 \frac{h}{h_0}\log_2\!\left(\frac{h_0}{h}\right),
 \qquad e=h\ebar
 \label{eq:criterion}
\end{equation}

where $S'$ is the integral approximation after bisection, and $h_0$ is the initial half-width.

The 9-point rule of AQNN9 has a fixed order, and its local error for a smooth function is roughly proportional to $h^{11}f^{(10)}$. Consequently, the number of function evaluations grows rapidly as the accuracy requirement is raised. In AQE11D, for any unconverged subinterval satisfying $|\ebar_9/\ebar_5|<P=0.2$, the selection parameter of~\cite{hibino2003} ($\ebar_9$ is the above $\ebar$ and $\ebar_5$ a 5-point-rule error functional built from the same nine values), nested 13-, 19-, 27-, and 41-point higher-order rules are applied in turn, which alleviates the lack of order. Their abscissae are constructed as symmetric binary fractions, and their weights are stored as rational numbers in MPFR variables. For problems whose endpoint value diverges, such as $1/\sqrt{x}$, we introduce a two-term treatment that never evaluates the endpoint itself, whose value in arbitrary precision is an infinity. For problems with a singularity at the endpoint, semi-analytic processing integrates the singular subinterval
analytically such that the local error can be kept below $\epsa$. In the original AQNN9, to keep the model error of the semi-analytic processing below $\epsa$, the singular subinterval must be refined down to $h\lesssim\epsa^{1/(p+2)}h_0$, where $p>-1$ is the exponent of the algebraic singularity.

The detection criteria are those of~\cite{ninomiya1980}; arbitrary precision requires only their thresholding and ordering to be changed. First, the rounding-error
noise floor of $\ebar$ is generalized from the fixed $32\varepsilon_{\mathrm{M}}$ of the double-precision version ($\varepsilon_{\mathrm{M}}$ being the machine epsilon of binary64) to a quantity proportional to $2^{5-P_B}$. Second, the order of the anomaly tests is reversed with respect to the double-precision version such that the arithmetic sequence test (logarithmic singularity) is performed before the constant sequence test (discontinuity). Because $\ebar_i$ forms an arithmetic sequence in the logarithmic singularity and the bisection depth required for detection increases with arbitrary precision (approximately 47 levels at $\epsa=10^{-35}$ against approximately 12 in double precision), the ratio $|\Delta\ebar|/|\ebar|$ decreases below the decision threshold. Thus, performing the constant-sequence test first would misclassify it as a
discontinuity. At a genuine discontinuity, $|\Delta\ebar|$ drops below the noise floor and therefore passes the arithmetic-sequence test. Thus, the reversed order correctly distinguishes between the two.

\subsection{Double Exponential Formula}
For a finite interval $[a,b]$, we apply a change of variable:

\begin{equation}
 x(t)=c+r\tanh\!\left(\frac{\pi}{2}\sinh t\right),
 \qquad
 c=\frac{a+b}{2},\quad r=\frac{b-a}{2}
 \label{eq:de-transform}
\end{equation}

The transformed integrand decays double exponentially as $|t|\to\infty$, so that the trapezoidal rule with step size $\Delta t$ and $2N+1$ sample points,

\begin{equation}
 \int_a^b f(x)\,dx
 \approx
 \Delta t\sum_{k=-N}^{N} f(x(k\Delta t))\,x'(k\Delta t)
 \label{eq:de-trapezoid}
\end{equation}

yields a highly accurate approximation ($\Delta t$ is unrelated to the $h$ of Sect.~2.1).

The endpoint treatment below is also new. In high-precision computation, one must avoid the cancellation or the evaluation of a divergent value resulting from $x(t)$ being rounded to the endpoint $a$ or $b$. We therefore compute the relative distance from
the endpoint,

\begin{equation}
 \delta=\frac{2}{1+\exp(2|u|)},
 \qquad u=\frac{\pi}{2}\sinh t,
\end{equation}
directly, and then generate the abscissae as $x=a+r\delta$ for $t<0$ and $x=b-r\delta$ for $t>0$.

While the DE formula is strong for functions with endpoint singularities and for smooth functions, it converges poorly on problems with an interior discontinuity, strong oscillation, or a narrow peak. Incorporating the algebraic treatment adopted in AQE11D into it is left for future work.

\section{Numerical Experiments}
The computing environment used in this study consists of the following GB10 machine:

\begin{description}
 \item[GB10] NVIDIA DGX Spark~\cite{nvidia_dgx_spark}, CPU (Cortex-X925 3.9\,GHz$\times 10$,
 Cortex-A725 2.81\,GHz$\times 10$), 128\,GB (LPDDR5X, unified memory) RAM, Ubuntu 24.04 LTS,
 GCC 13.3.0, MPFR 4.3.0, GMP 6.4.1
\end{description}

All runs are performed on a CPU with OpenMP. We set the accuracy requirements to $\epsa=10^{-50}$ and $10^{-100}$, with mantissa precisions $P_B=256$ and $512$ bits, respectively. The safety limit on the number of function evaluations is $N_{\max}=10^{8}$ and the maximum bisection depth is 500. The tolerance test employs the absolute error $\epsa$ specified in advance by the user. Kahaner's 21 test problems used in the experiments are listed in \tablename~\ref{tab:probs}.

\begin{table}[tb]\centering\footnotesize
\setlength{\tabcolsep}{7pt}
\caption{Kahaner's 21 test problems~\cite{kahaner1971}. $I$: analytic where an elementary antiderivative
exists, else mpmath at 240 digits, rounded to 7 decimals; we use the exact $\pi$, unlike~\cite{ninomiya1980,
hibino2003}.  Problems 5, 20: analytic, nearby complex poles (sing.\ = singularity, osc.\ = oscillation).}
\label{tab:probs}
\resizebox{\columnwidth}{!}{%
\begin{tabular}{@{}r @{\ }l l r l@{}}\toprule
\# & $f(x)$ & interval & exact value $I$ & character \\ \midrule
1 & $e^x$ & [0,1] & $1.7182818$ & smooth \\
2 & $1\,(x>0.3),\,0$ & [0,1] & $0.7000000$ & discontinuity \\
3 & $\sqrt{x}$ & [0,1] & $0.6666667$ & alg.\ sing.\ (weak) \\
4 & $\tfrac{23}{25}\cosh x-\cos x$ & [-1,1] & $0.4794282$ & smooth \\
5 & $1/(x^4+x^2+0.9)$ & [-1,1] & $1.5822330$ & peak, cplx.\ poles \\
6 & $x^{3/2}$ & [0,1] & $0.4000000$ & alg.\ sing.\ (weak) \\
7 & $1/\sqrt{x}$ & [0,1] & $2.0000000$ & alg.\ sing.\ (strong) \\
8 & $1/(1+x^4)$ & [0,1] & $0.8669730$ & smooth \\
9 & $2/(2+\sin 10\pi x)$ & [0,1] & $1.1547005$ & oscillatory \\
10 & $1/(1+x)$ & [0,1] & $0.6931472$ & smooth \\
11 & $1/(1+e^x)$ & [0,1] & $0.3798855$ & smooth \\
12 & $x/(e^x-1)$ & [0,1] & $0.7775046$ & removable sing. \\
13 & $\sin(100\pi x)/(\pi x)$ & [0.1,1] & $0.0090986$ & strong osc. \\
14 & $\sqrt{50}\,e^{-50\pi x^2}$ & [0,10] & $0.5000000$ & sharp peak \\
15 & $25 e^{-25x}$ & [0,10] & $1.0000000$ & sharp peak \\
16 & $50/(\pi(2500x^2+1))$ & [0,10] & $0.4993634$ & sharp peak \\
17 & $50(\sin 50\pi x/50\pi x)^2$ & [0.01,1] & $0.1121393$ & sinc$^2$ osc. \\
18 & \begin{tabular}{@{}l@{}}$\cos(\cos x+3\sin x+2\cos 2x$\\$\qquad{}+3\sin 2x+3\cos 3x)$\end{tabular} & [0,$\pi$] & $0.8386763$ & strong osc. \\
19 & $\ln x$ & [0,1] & $-1.0000000$ & log sing. \\
20 & $1/(1.005+x^2)$ & [-1,1] & $1.5643964$ & cplx.\ poles \\
21 & $\sum_{i=1}^{3}\operatorname{sech}^{2i}\!\left(10^{i}(x-\tfrac{i}{5})\right)$ & [0,1] & $0.2108027$ & 3 sharp peaks \\
\bottomrule\end{tabular}}\end{table}

\subsection{Attained Accuracy and Function Evaluations}
The attained accuracy and number of function evaluations are shown in \tablename~\ref{tab:acc}. Over the range we compared by raising $P_B$ from 256 and 512 up to 1024 bits, the convergence results for the same tolerance are almost identical, such that, under the present conditions, the order of the quadrature rule and the stopping criterion dominate the working precision. As this holds only when the rounding error is sufficiently
small, we interpret it as a trend within the experimental conditions rather than as genuine independence of the working precision.

AQE11D reaches the target accuracy on all 21 problems at both tolerances. In particular, for problem 7, $1/\sqrt{x}$, the treatment that avoids evaluating the endpoint produces an error of $3.8\times10^{-75}$ at
$\epsa=10^{-50}$. However, at $\epsa=10^{-100}$, when the algebraic singularity at the left endpoint $x=0$ (exponent $-1/2$) is detected and treated with two-term processing, the singular subinterval must be refined down to $h\lesssim\epsa^{1/(p+2)}h_0$ (i.e., \ $\epsa^{2/3}h_0$ for $p=-1/2$, a bisection depth of about 220 levels) to keep the model
error of the semi-analytic processing below $\epsa$. The function evaluations on the sequence of subintervals leading there reach about $5.4\times10^{7}$ (yielding an error of $6.6\times10^{-150}$). This illustrates the cost limit of a polynomial-type adaptive rule handling a strong endpoint singularity at high accuracy: The DE formula solves the same problem with 743 evaluations, reaching $2.4\times10^{-105}$, a cost difference of about $7\times10^{4}$. (With $N_{\max}=2\times10^{7}$, problem 7 is cut off at $1.8\times10^{-10}$, and we therefore use $10^{8}$ here.) 

The DE formula converges on 18 problems at both settings; it does not converge on problem 2 (interior discontinuity), problem 9 (strong oscillation), and problem 21 (sharp peaks). Overall, on the problems where it converges, the DE formula needs fewer function evaluations and is faster.

\begin{table}[tb]\centering\footnotesize
\setlength{\tabcolsep}{7pt}
\caption{Attained accuracy $|$err$|$ and number of function evaluations Nfe (problem numbers follow
Table~\ref{tab:probs}).  $^{*}$: DE did not converge.  $^{\dagger}$: strong endpoint singularity,
very costly (see text).}
\label{tab:acc}
\resizebox{\columnwidth}{!}{%
\begin{tabular}{@{}r rr rr rr rr@{}}\toprule
 & \multicolumn{4}{c}{$\varepsilon_a=10^{-50}$ ($P_B=256$)} & \multicolumn{4}{c}{$\varepsilon_a=10^{-100}$ ($P_B=512$)}\\
\cmidrule(lr){2-5}\cmidrule(lr){6-9}
 & \multicolumn{2}{c}{AQE11D} & \multicolumn{2}{c}{DE} & \multicolumn{2}{c}{AQE11D} & \multicolumn{2}{c}{DE}\\
\cmidrule(lr){2-3}\cmidrule(lr){4-5}\cmidrule(lr){6-7}\cmidrule(lr){8-9}
\# & $|$err$|$ & Nfe & $|$err$|$ & Nfe & $|$err$|$ & Nfe & $|$err$|$ & Nfe \\ \midrule
1 & $1.03e-77$ & 81 & $1.86e-55$ & 289 & $1.07e-154$ & 8901 & $1.13e-105$ & 655 \\
2 & $9.28e-53$ & 1681 & $5.95e-02$$^{*}$ & 95033 & $9.92e-103$ & 3341 & $5.95e-02$$^{*}$ & 100441 \\
3 & $1.73e-74$ & 29641 & $5.01e-56$ & 289 & $5.32e-149$ & 2707751 & $3.03e-106$ & 655 \\
4 & $1.92e-78$ & 221 & $3.13e-53$ & 567 & $3.51e-154$ & 16741 & $1.07e-105$ & 655 \\
5 & $1.12e-77$ & 2125 & $7.64e-52$ & 1117 & $1.37e-153$ & 143245 & $2.25e-101$ & 2573 \\
6 & $3.44e-73$ & 13815 & $5.01e-56$ & 289 & $7.15e-148$ & 1019337 & $3.03e-106$ & 655 \\
7 & $3.83e-75$ & 182181 & $2.03e-55$ & 333 & $6.61e-150$$^{\dagger}$ & 54240551 & $2.37e-105$ & 743 \\
8 & $6.22e-74$ & 1691 & $1.34e-53$ & 567 & $2.95e-150$ & 100811 & $6.55e-104$ & 1297 \\
9 & $3.82e-75$ & 21323 & $1.35e-01$$^{*}$ & 95045 & $9.79e-149$ & 1283423 & $1.35e-01$$^{*}$ & 100449 \\
10 & $8.28e-76$ & 851 & $1.34e-53$ & 567 & $3.34e-150$ & 54191 & $6.55e-104$ & 1297 \\
11 & $8.86e-79$ & 501 & $6.85e-54$ & 567 & $1.24e-150$ & 17861 & $3.36e-104$ & 1297 \\
12 & $6.44e-78$ & 221 & $7.93e-56$ & 289 & $3.82e-150$ & 12611 & $4.80e-106$ & 655 \\
13 & $8.12e-78$ & 34883 & $2.62e-53$ & 1907 & $7.28e-155$ & 2279713 & $1.47e-103$ & 2239 \\
14 & $3.40e-55$ & 3817 & $2.73e-51$ & 2239 & $4.00e-108$ & 238301 & $3.38e-101$ & 5147 \\
15 & $7.58e-55$ & 3867 & $5.11e-52$ & 1129 & $9.71e-106$ & 240897 & $2.57e-103$ & 1301 \\
16 & $2.66e-75$ & 9121 & $2.50e-51$ & 2241 & $1.48e-149$ & 557431 & $2.91e-101$ & 5153 \\
17 & $2.13e-78$ & 34547 & $3.15e-52$ & 2241 & $4.06e-154$ & 2131817 & $1.29e-102$ & 2583 \\
18 & $6.30e-78$ & 7653 & $9.46e-52$ & 1119 & $3.41e-152$ & 436963 & $4.33e-102$ & 2577 \\
19 & $2.52e-75$ & 62451 & $1.81e-57$ & 291 & $1.15e-149$ & 7564241 & $3.44e-107$ & 657 \\
20 & $2.21e-77$ & 2125 & $1.10e-51$ & 1119 & $1.15e-153$ & 143245 & $3.25e-101$ & 2573 \\
21 & $3.44e-61$ & 19781 & $1.38e-02$$^{*}$ & 94923 & $2.67e-127$ & 1209963 & $1.38e-02$$^{*}$ & 100419 \\
\bottomrule\end{tabular}}\end{table}

\subsection{Serial Execution Time per Problem}
\begin{figure}[tb]
 \centering
 \includegraphics[width=\linewidth]{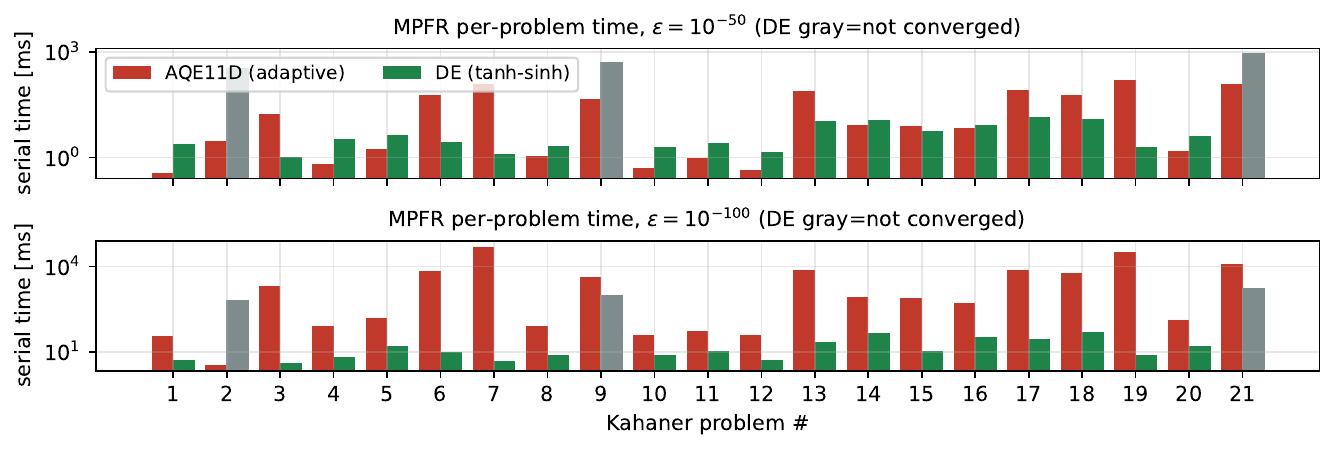}
 \caption{Single-core execution time per problem (MPFR; logarithmic vertical axis). Upper:
 $\epsa=10^{-50}$ ($P_B=256$); lower: $\epsa=10^{-100}$ ($P_B=512$). AQE11D is shown in red, the converged DE formula in green, and the unconverged DE formula in gray.}
 \label{fig:time}
\end{figure}

\figurename~\ref{fig:time} shows the computation time of each problem executed on one core. On the endpoint-singular problems 3, 7 and 19 and the strongly oscillatory problems 13 and 17, the DE formula converges in a far shorter time than AQE11D does: At $\epsa=10^{-50}$, problem 7 takes 124.4~ms with AQE11D and 1.2~ms with the DE formula, and problem 19 takes 158.4~ms and 2.0~ms, respectively. By contrast, on problems 2, 9, and 21, as mentioned earlier, the DE formula does not converge and AQE11D is preferable.

\subsection{Parallel Performance}
We evaluate parallelization at two granularities: intra-integral parallelism (fine-grained), in which the function evaluations within each integral are parallelized with OpenMP, and problem-wise task parallelism
(coarse-grained), in which independent integration problems are assigned to the threads. In the intra-integral case, AQE11D evaluates the six newly required points in parallel at each bisection, while the DE formula evaluates the sample points of each level in parallel in blocks. In both cases, the sums are accumulated in sequential order; thus, exploiting an MPFR property, the results are identical to the serial ones.

\begin{figure}[tb]
 \centering
 \includegraphics[width=\linewidth]{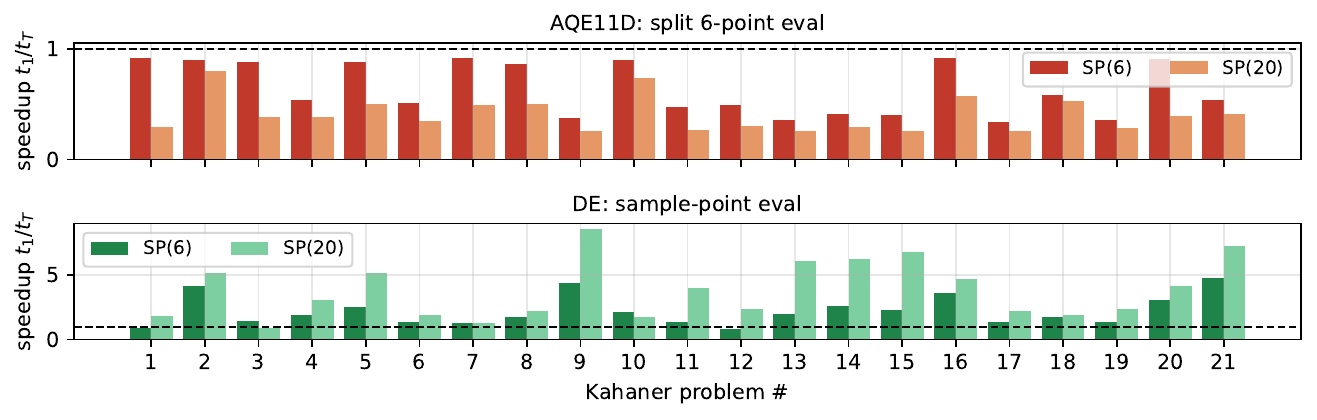}
 \caption{Speedup of the parallelized function evaluation per problem ($\epsa=10^{-100}$, $P_B=512$).
 Upper: AQE11D (six points per bisection); lower: DE (sample points). The dashed line is $\mathrm{SP}=1$.}
 \label{fig:inner}
\end{figure}

\figurename~\ref{fig:inner} shows the speedup $\mathrm{SP}(T)=t_1/t_T$ of the intra-integral parallelism when each problem is run on its own. For AQE11D, the degree of parallelism is limited to six points, and the overhead of creating and synchronizing a \texttt{parallel for} at every bisection, together with that of creating and destroying MPFR variables, is relatively large.

Across all problems, $\mathrm{SP}(6)\le0.92$ (no speedup at all), and 20 threads make it worse still. The DE formula, by contrast, shows higher parallelism the more sample points a problem has per level: $\mathrm{SP}(6)\approx4$ and $\mathrm{SP}(20)\approx5$--$9$ on the unconverged problems (2, 9, 21), the sample points of which reach about $10^{5}$, and $\mathrm{SP}(6)\approx2$--$3.6$ on the peaked problems (5, 14, 15, 16 and 20).

On the converged problems, with only a few hundred sample points, even the DE formula shows a modest $\mathrm{SP}$. In other words, the uniform sample-point evaluation of the DE formula suits fine-grained parallelism, whereas the adaptive six-point bisection of AQE11D gains little from it. The largest speedup per method and problem is also listed in \tablename~\ref{tab:time}.

\tablename~\ref{tab:time} presents the total execution time of the 21 problems under problem-wise task parallelism. At $\epsa=10^{-50}$, AQE11D obtains a speedup of about 4.8 on 20 cores, but at $\epsa=10^{-100}$, this decreases to 2.6; this is because problem 7 alone, with its strong endpoint singularity, takes about 50~s serially and, even under dynamic scheduling, is assigned to a single core and thus determines the total time of the
20-core run (Amdahl's law). For the DE formula, the three unconverged problems (2, 9, 21) cannot be subdivided; thus, the speedup stays around 2 at both precisions. Since the DE formula solves the singular problems with a few hundred evaluations, its total execution time is short, and its
advantage over AQE11D widens as the precision is raised, reaching a factor of about 40 at $\epsa=10^{-100}$.
\IfFileExists{tab_time_en.tex}{\begin{table}[tb]\centering\small
\caption{Total time [s] for the 21-problem battery (20-core CPU, task parallelism, $N_{\max}=10^{8}$).
``Speedup'' is the 20-core speedup of the task parallelism; ``max SP'' the largest speedup of the
intra-integral parallelism (Fig.~\ref{fig:inner}).}
\label{tab:time}
\resizebox{\columnwidth}{!}{%
\begin{tabular}{@{}ll rr r r@{}}\toprule
$\varepsilon_a$ & method & serial [s] & 20 cores [s] & speedup & max SP \\ \midrule
\multirow{2}{*}{$10^{-50}$} & AQE11D & 0.77 & 0.16 & 4.84$\times$ & 0.97$\times$ \\
 & DE & 1.77 & 0.88 & 2.02$\times$ & 6.74$\times$ \\ \midrule
\multirow{2}{*}{$10^{-100}$} & AQE11D & 129.7 & 50.4 & 2.57$\times$ & 0.92$\times$ \\
 & DE & 3.25 & 1.55 & 2.09$\times$ & 8.54$\times$ \\
\bottomrule\end{tabular}}\end{table}
}{}

\section{Conclusion and Future Work}
We implemented arbitrary-precision AQE11D and the DE formula on a CPU using MPFR and evaluated their accuracy and parallel performance on Kahaner's 21 test problems. AQE11D reaches the target on all 21 problems at both precisions but needs a large number of evaluations for a strong endpoint singularity. The DE formula is strong against endpoint singularities and converges on 18 problems at both precisions; on these problems, it needs far fewer function evaluations and far less time than AQE11D.

In parallelization, task parallelism gives AQE11D a speedup of around 4.8 (only 2.6 at $\epsa=10^{-100}$, where the strongly endpoint-singular problem 7 presents a bottleneck), and parallelizing the function evaluations within an integral is effective only for the sample-point evaluation of the DE formula. 

In future work, we intend to address parallelization at the level of subintervals, better load balancing, hybridization at the interval level, and throughput improvement using GPU batch processing with \texttt{mpc\_cuda}~\cite{mpccuda}, a CUDA port of MPFR/MPC~\cite{mpc}, aiming to develop a numerical integrator applicable to a wider range of problems.

\section*{Acknowledgments}
We thank the anonymous reviewer for their helpful comments on the original manuscript. This work was supported by JSPS KAKENHI Grant Number 26K14846.

\end{document}